\documentstyle[amstex,12pt]{article}
\title{Generalized Umemura polynomials}
\author{Anatol N. Kirillov and Makoto Taneda}

\date{}

\begin{document}

\newcommand{\bbb}{{\bf b}}%
\newcommand{\C}{{\mathbb C}}%
\newcommand{\HH}{{\cal H}}%
\newcommand{\dd}{{\bf d}}
\newcommand{\cc}{{\bf c}}
\newcommand{\odd}{{\rm odd}}
\newcounter{newsection}
\def\neweq{\setcounter{equation}{0}}
\renewcommand{\thesection}{\S \arabic{section}.}
\renewcommand{\thesubsection}{\thesection.\arabic{subsection}}
\renewcommand{\theequation}{\thenewsection.\arabic{equation}}
\renewcommand
{\thenewsection}{\setcounter{newsection}{\value{section}}\arabic{newsection}}
\newtheorem{theorem}{Theorem}[section]
\newtheorem{Proposition}[theorem]{Proposition}
\newtheorem{Definition}[theorem]{Definition}
\newtheorem{Theorem}[theorem]{Theorem}
\newtheorem{Remarks}[theorem]{Remarks}
\newtheorem{Corollary}[theorem]{Corollary}
\newtheorem{Conjecture}[theorem]{Conjecture}
\newtheorem{Lemma}[theorem]{Lemma}
\maketitle
\begin{abstract}
We introduce and study generalized Umemura polynomials $U_{n,m}^{(k)}(z,w;a,b)$
which are a natural generalization of the Umemura polynomials $U_n(z,w;a,b)$
related to Painlev\'e $VI$ equation.
We will show that if $a=b$, or $a=0$, or $b=0$ then polynomials
$U_{n,m}^{(0)}(z,w;a,b)$ generate solutions to Painlev\'e $VI$.
We will describe a connection between polynomials $U_{n,m}^{(0)}(z,w;a,0)$
and certain Umemura polynomials $U_k(z,w;\alpha,\beta)$.
\end{abstract}
\section{Introduction}
There exists a vast body of literature about the Painlev\'e $VI$ equation
$P_{VI}:=P_{VI}(\alpha,\beta,\gamma,\delta):$
\begin{eqnarray}
\frac{d^2q}{dt^2}&=&\frac{1}{2}
\left(\frac{1}{q}+\frac{1}{q-1}+\frac{1}{q-t}\right)
\left(\frac{dq}{dt}\right)^2
-
\left(\frac{1}{t}+\frac{1}{t-1}+\frac{1}{q-t}\right)
\left(\frac{dq}{dt}\right) \nonumber \\
\smallskip \nonumber\\
&&+
\frac{q(q-1)(q-t)}{t^2(t-1)^2}
\left(
\alpha-\beta\frac{t}{q^2}
+\gamma\frac{(t-1)}{(q-1)^2}
+\delta\frac{t(t-1)}{(q-t)^2}
\right)
\end{eqnarray}
where $t \in {\bf C}$, $q:=q(t;\alpha,\beta,\gamma,\delta)$ is
a function of $t$, and $\alpha,\beta,\gamma,\delta$ are arbitrary complex
parameters.
It is well-known and goes back to Painlev\'e that any solution $q(t)$ of
the equation $P_{VI}$ satisfies the so-called Painlev\'e property:
\begin{itemize}
\item
the critical points $0,1$ and $\infty$ of the equation (1,1) are
the only \emph{fixed singularities} of $q(t)$.
\item
any \emph{movable singularity} of $q(t)$
(the position of which depends on integration constants)
is a pole.
\end{itemize}
In this paper we introduce and initiate the study of certain
special polynomials related to the Painlev\'e $VI$ equation,
namely, the generalized Umemura polynomials
$U_{n,m}^{(k)}(z,w;a,b)$. These polynomials have many interesting
combinatorial and algebraic properties and in the particular case
$n=0=k$ coincide with Umemura's polynomials $U_m(z^2,w^2;a,b)$,
see e.g. [U], [NOOU]. In the present paper we study recurrence
relations between polynomials $U_{n,m}^{(k)}(z,w;a,b)$. Our main
result is Theorem~\ref{t4.1} which gives a generalization of the
recurrence relation between Umemura's polynomials [U]. In some
particular cases the recurrence relation obtained in
Theorem~\ref{t4.1} coincides with that for Umemura's polynomials.
As a corollary, we obtain that polynomials
$U_{n,m}^{(0)}(z,w;a,0)$ also generate solutions to the equation
Painlev\'e $VI$. The main mean in our proofs is Lemma~\ref{l4.2}
from Section 4. For example, we prove a new recurrence relation
between Umemura's polynomials (Theorem~\ref{t4.9}), and describe
explicitly connections between polynomials $U_{n,m}(0,b)$ and
Umemura's polynomials $U_m(b_1,b_2)$, see Lemma~\ref{l4.7}.
Finally, in section 5 we state a conjectural description of the
Pl\"ucker relations between certain Umemura's polynomials.

\section{Painlev\'e $VI$}
\setcounter{equation}{0}
\renewcommand{\thesection}{\arabic{section}}
In this section
we collect together some basic results about equation Painlev\'e $VI$.
More detail and proofs may be found in familiar series of papers by
K.Okamoto [OI-OIV].
We refer the reader to the Proceedings of Conference "The Painlev\'e property.
One century later" [P], where different aspects of the theory of Painlev\'e
equations may be found.

\subsection{Hamiltonian form}
It is well-known and goes back to a paper by Okamoto [OI] that
the sixth Painlev\'e equation (1.1) is equivalent to the
following Hamiltonian system
\begin{eqnarray}
\HH_{VI}( \bbb ;t,q,p):
\left\{
\begin{array}{lcr}
\frac{\displaystyle dq}{\displaystyle dt} &=&
\frac{\displaystyle \partial H}{\displaystyle \partial p},\\
\smallskip \\
\frac{\displaystyle dp}{\displaystyle dt} &=&
-\frac{\displaystyle \partial H}{\displaystyle \partial q}
\end{array}
\right.
\end{eqnarray}
with the Hamiltonian
\begin{eqnarray}
\lefteqn{H:=H_{VI}(\bbb;t,q,p)
=\frac{1}{t(t-1)}
\left[
q(q-1)(q-t)p^2
\right.
- \left\{
(b_1+b_2)(q-1)(q-t)
\right.
}
\nonumber\\
\smallskip \nonumber\\
&&
\left. \left.
+(b_1-b_2)q(q-t)
+(b_3+b_4)q(q-1)
\right\}p
+(b_1+b_3)(b_1+b_4)(q-t)
\right],\nonumber
\end{eqnarray}
where $\bbb=(b_1,b_2,b_3,b_4)$ belongs to the parameters space $\C^4$;
the parameters $(\alpha,\beta, \gamma, \delta)$ and $(b_1,b_2,b_3,b_4)$
are connected by the following relations
\begin{eqnarray}
\alpha=\frac{1}{2}(b_3-b_4)^2,
\beta=-\frac{1}{2}(b_1+b_2)^2,
\gamma=\frac{1}{2}(b_1-b_2)^2,
\delta=-\frac{1}{2}(b_3-b_4)(b_3+b_4-2).
\end{eqnarray}
\begin{Proposition}\label{p2.1}(K. Okamoto, [OI].)
If $(q(t),p(t))$ is a solution to the Hamiltonian system (2.1),
the function
\begin{eqnarray*}
h(\bbb,t)=t(t-1)H_{VI}(\bbb;t,q(t),p(t))
+e_2(b_1,b_3,b_4)-\frac{1}{2}e_2(b_1,b_2,b_3,b_4)
\end{eqnarray*}
satisfies the equation $E_{VI}(\bbb):$
\begin{eqnarray}
\frac{dh}{dt} \left[t(t-1)\frac{d^2h}{dt^2}\right]^2
+\left[\frac{dh}{dt}\left\{2h-(2t-1)\frac{dh}{dt}
\right\}+b_1b_2b_3b_4\right]
=\prod_{k=1}^{4}\left(\frac{dh}{dt}+b_k^2\right).
\end{eqnarray}
Conversely, for a solution $h:=h(\bbb,t)$ to the equation $E_{VI}(\bbb)$
such that $\frac{d^2h}{dt^2} \neq 0$, there exists a solution $(q(t),p(t))$ to
the Hamiltonian system (2.1).
Furthermore, the function $q:=q(t)$ is a solution to the Painlev\'e equation
(1.1), where parameters $(\alpha, \beta, \gamma, \delta)$ are determined
by the relations (2.2).

\end{Proposition}
We will call the equation $E_{VI}(\bbb)$ by the
Painlev\'e--Okamoto equation.

\subsection{B\"acklund transformation}
Consider the following linear transformation of the parameters space $\C^4$:
\begin{eqnarray*}
&&s_1:=(b_1,b_2,b_3,b_4)\longmapsto (b_2,b_1,b_3,b_4),\\
\smallskip \\
&&s_2:=(b_1,b_2,b_3,b_4)\longmapsto (b_1,b_3,b_2,b_4),\\
\smallskip \\
&&s_3:=(b_1,b_2,b_3,b_4)\longmapsto (b_1,b_2,b_4,b_3),\\
\smallskip \\
&&s_0:=(b_1,b_2,b_3,b_4)\longmapsto (b_1,b_2,-b_3,-b_4),\\
\smallskip \\
&&l_3:=(b_1,b_2,b_3,b_4)\longmapsto (b_1,b_2,b_3+1,b_4).
\end{eqnarray*}
Denote by $W=<s_0,s_1,s_2,s_3,l_3>$ the subgroup of ${\rm Aut} \C^4$
generated by these transformation.
It is not difficult to see, that $W \cong W(D_4^{(1)})$, i.e.
$W$ is isomorphic to the affine Weyl group of type $D_4^{(1)}$.
\begin{Proposition}\label{p2.2}(K. Okamoto, [OI].)
For each $w \in W$, there exists a birational transformation
\begin{eqnarray*}
L_w:\{\mbox{solutions to ${\cal H}_{VI}(\bbb)$ } \} \longmapsto
\{\mbox{solutions to ${\cal H}_{VI}(w(\bbb))$ } \}.
\end{eqnarray*}
\end{Proposition}
The birational transformations $L_w$, $w \in W(D_4^{(1)})$
are called by {\it B\"acklund transformations} associated to the equation
Painlev\'e $VI$.

\subsection{$\tau$-function}
Let $(q(t),p(t))$ be a solution to the Hamiltonian system (2.1),
the $\tau$-function $\tau(t)$ corresponding to the solution $(q(t),p(t))$
is defined by the following equation
\begin{eqnarray*}
\frac{d}{dt} \log \tau(t) =H_{VI}(\bbb;t,q(t),p(t));
\end{eqnarray*}
in other words,
\begin{eqnarray*}
\tau(t)=({\rm constant}) \exp \left( \int H_{VI}(\bbb;t,q(t),p(t)) dt\right).
\end{eqnarray*}

\subsection{Umemura polynomials}
Suppose that $b_3=-\frac{1}{2}$, $b_4=0$, then
it is well-known and goes back to Umemura's paper [U]
, that the pair
\begin{eqnarray*}
(q_0,p_0)=\left(
\frac{(b_1+b_2)^2-(b_1^2-b_2^2)\sqrt{t(1-t)}}{(b_1-b_2)^2+4 b_1 b_2 t},
\frac{b_1q_0-\frac{1}{2}(b_1+b_2)}{q_0(q_0-1)}
\right)
\end{eqnarray*}
defines a solution to the Hamiltonian system (2.1) with parameters
$\bbb=(b_1,b_2,-\frac{1}{2},0)$.
Note, see e.g. [U], that
\begin{eqnarray*}
\lefteqn{H_0=H_{VI}\left((b_1,b_2,-\frac{1}{2},0);t,q_0(t),p_0(t)\right)}\\
\smallskip \\
&=&\frac{1}{t(t-1)}
\left\{
b_1(b_1-1)(1-2t)+2 b_1^2 \sqrt{t(t-1)}+2b_2(b_1-t) \right.\\
\smallskip \\
&&\left.
 +b_2(b_2-1)(1-2t)
-2b_2^2 \sqrt{t(t-1)}
\right\},
\end{eqnarray*}
and \begin{eqnarray*}
\tau_0(t)=\exp
\left\{\int H_0(t)dt
\right\}.
\end{eqnarray*}
To introduce Umemura's polynomials,
let $(q_m,p_m)$ be a solution to the Hamiltonian system
${\cal H}_{VI}(b_1,b_2,-\frac{1}{2}+m,0)=
{\cal H}_{VI}(l_3^m(b_1,b_2,-\frac{1}{2},0)$
obtained from the solution $(q_0,p_0)$ by applying $m$ times the
the B\"acklund transformation $l_3$.
Consider the corresponding $\tau$-function $\tau_m$:
\begin{eqnarray*}
\frac{d}{dt} \log \tau_m =H_{VI}((b_1,b_2,-\frac{1}{2}+m,0);t,q_m(t),p_m(t)).
\end{eqnarray*}
It follows from Proposition~\ref{p2.1}, see e.g. [OI], [U], that
the $\tau$-functions $\tau_n:=\tau_n(t)$ satisfy the Toda equation
\begin{eqnarray}
\frac{\tau_{n-1} \tau_{n+1}}{\tau_n^2}=
\frac{d}{dt}
\left( t(t-1)\frac{d}{dt}(\log \tau_n)
\right)+(b_1+b_2+n)(b_3+b_4+n).
\end{eqnarray}
Follow H. Umemura [U], define a family of functions $T_n(t)$, $n=0,1,2,\ldots$,
by
\begin{eqnarray*}
\tau_n(t)=T_n(t) \exp
\left( \int
\left( H_0(t)-\frac{n(b_1t-\frac{1}{2}(b_1+b_2))}{t(t-1)} \right)dt
\right).
\end{eqnarray*}
\begin{Proposition}\label{p2.3}(H. Umemura, [U])
$T_n(t)$ is a polynomial in variable
$v:=\sqrt{\frac{t}{t-1}}+\sqrt{\frac{t-1}{t}}$
with rational coefficients.
\end{Proposition}
For example, $T_0=1$, $T_1=1$,
$T_2=\frac{1}{2}\left( -4b_1^2+1)(2-v)/4+(-4b_2^2+1)(2+v)/4\right)$.\\
It follows from the Toda equation (2.4) that polynomials $T_n:=T_n(v)$
satisfy the following recurrence relation [U]:
\begin{eqnarray}
T_{n-1}T_{n+1} &=& \left \{
\frac{1}{4}(-2b_1^2-2b_2^2+(b_1^2-b_2^2)v)+(n-\frac{1}{2})^2
\right\}
T_n^2 \\
\smallskip \nonumber\\
&& + \frac{1}{4}(v^2-4)^2\left\{T_n \frac{d^2T_n}{dv^2}
-\left(\frac{d T_n}{dv}\right)^2\right\} \nonumber \\
\smallskip \nonumber\\
&& + \frac{1}{4}(v^2-4)v T_n \frac{d T_n}{dv}\nonumber
\end{eqnarray}
with initial conditions $T_0=T_1=1$.
\begin{Definition}\label{d2.4}
Polynomials $U_n:=U_n(z,w,b_1,b_2):=2^{n(n-1)}T_n(v)$,
where $z=\frac{2-v}{4}$, $w=\frac{2+v}{4}$, are called by Umemura polynomials.
\end{Definition}
The formula (2.6) below was stated as a conjecture by
S. Okada, M. Noumi, K. Okamoto and H. Umemura [NOOU] and has been proved
recently by M. Taneda, and A. N. Kirillov (independently):
\begin{eqnarray}
2^{n(n-1)}T_n(v) :=U_n(z,w,b_1,b_2)= \sum_{I \subset[n-1]}
\dd_n(I)
c_I d_{[n-1] \backslash I}
 \ z^{|I|} \  w^{|I^{c}|},
\end{eqnarray}
where
\begin{description}
\item[(i)] $[n-1]=\{1,2,\ldots,n-1\}$;
for any subset $I=\{i_1>i_2>\cdots>i_p\} \subset [n-1]$,
$\dd_n(I)=\dim_{\lambda(I)}^{GL(n)}$ stands for
the dimension of irreducible representation of
the general linear group $GL(n)$ corresponding to
the highest weight $\lambda(I)$
with the Frobenius' symbol
$\lambda(I)=(i_1,i_2, \ldots,i_p|i_1-1,i_2-1, \ldots, i_p-1)$;
\item[(ii)]$c=-4b_1^2$, $d=-4b_2^2$, $z=\frac{2-v}{4}$, $w=\frac{2+v}{4}$;
\item[(iii)]$\bar{c}=c+(2k-1)^2$, $\bar{d}=d+(2k-1)^2$,
$c_k=\bar{c}_1 \bar{c}_2 \cdots \bar{c}_k$,
$d_k=\bar{d}_1 \bar{d}_2 \cdots \bar{d}_k$;
\item[(iv)]
$|I|=i_1+i_2+\cdots+i_p$.
\end{description}
Recall that Frobenius' symbol $(a_1,a_2,\ldots,a_p|b_1,b_2,\ldots,b_p)$
denotes the partition which corresponds to the following diagram
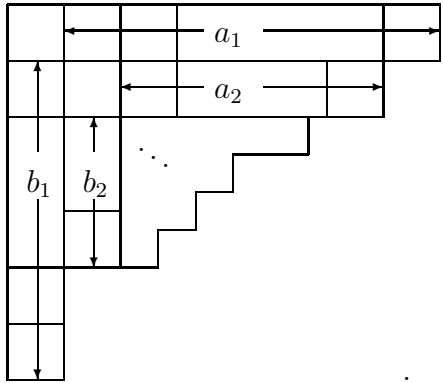
\begin{center}
 \vskip3mm
 \setlength{\unitlength}{0.5mm}
\begin{picture}(120,100)
\put(10,100){\line(1,0){115}}
\put(10,85){\line(1,0){115}}
\put(10,70){\line(1,0){100}}
\put(25,45){\line(1,0){15}}
\put(10,30){\line(1,0){40}}
\put(10,15){\line(1,0){15}}
\put(10,0){\line(1,0){15}}

\put(50,30){\line(0,1){10}}
\put(50,40){\line(1,0){10}}
\put(60,40){\line(0,1){10}}
\put(60,50){\line(1,0){10}}
\put(70,50){\line(0,1){10}}
\put(70,60){\line(1,0){20}}
\put(90,60){\line(0,1){10}}

\put(10,100){\line(0,-1){100}}
\put(25,100){\line(0,-1){100}}
\put(40,100){\line(0,-1){70}}
\put(55,100){\line(0,-1){30}}
\put(110,100){\line(0,-1){15}}
\put(125,100){\line(0,-1){15}}
\put(95,85){\line(0,-1){15}}
\put(110,85){\line(0,-1){15}}

\put(62,93){\vector(-1,0){37}}
\put(78,93){\vector(1,0){47}}
\put(65,90){$a_1$}

\put(62,78){\vector(-1,0){22}}
\put(78,78){\vector(1,0){32}}
\put(65,75){$a_2$}

\put(18,61){\vector(0,1){24}}
\put(18,48){\vector(0,-1){48}}
\put(15,50){$b_1$}

\put(33,61){\vector(0,1){9}}
\put(33,48){\vector(0,-1){18}}
\put(30,50){$b_2$}

\put(44,56){$\ddots$}
\put(115,0){.}

\end{picture}
\end{center}

\section{Generalized Umemura polynomials}
\setcounter{equation}{0}
\renewcommand{\thesection}{\arabic{section}}
Let $n,m,k$ be fixed nonnegative integers, $k \leq n$.
Denote by $[n;m]$ the set of integers $\{1,2,\ldots,n,n+2,n+4,\ldots,n+2m \}$.
Let $I$ be a subset of the set $[n;m]$. Follow [DK], define the numbers
\begin{eqnarray}
\dd_{n,m}(I)=\prod_{i \in I, j \in [n;m]\backslash I}
\left|\frac{i+j}{i-j}\right|,
\ \ \ \ \ \cc(I)=\sum_{i \in I, i>n} \frac{i-n}{2}
\end{eqnarray}
It has been shown in [DK], that in fact $\dd_{n,m}(I)$ are integers
for any subset $I \subset [n;m]$.
Now we are going to introduce the generalized Umemura polynomials
\begin{eqnarray*}
U_{n,m}^{(k)}:=U_{n,m}^{(k)}(z,w;a,b)=
\sum_{[k] \subset I \subset [n;m]}
\ \
\prod_{i \in I \backslash [k], j \in [k]}
\left(\frac{i+j}{i-j}\right)
\dd_{n,m}(I)
(-1)^{\cc(I)} e_{I}^{(n,m,k)}(z,w),
\end{eqnarray*}
where \\
\begin{description}
\item[(i)] $[k]$ stands for the set $\{1,2,\ldots,k\}$;
\item[(ii)]
$\bar{a}_k=a+(k-1)^2$, $\bar{b}_k=b+(k-1)^2$ and
$a_{2k}=\bar{a}_2 \bar{a}_4 \cdots \bar{a}_{2k}$,
$a_{2k+1}=\bar{a}_1 \bar{a}_3 \cdots \bar{a}_{2k+1}$;
$b_{2k}=\bar{b}_2 \bar{b}_4 \cdots \bar{b}_{2k}$,
$b_{2k+1}=\bar{b}_1 \bar{b}_3 \cdots \bar{b}_{2k+1}$;
\item[(iii)] for any subset $I \subset [n;m]$, we set
$a_I =\prod_{i \in I} a_i $, $b_I = \prod_{i \in I} b_i$;
\item[(iv)] $e_I^{(n,m,k)}(z,w)=a_{I \backslash [k]}
b_{[n,m]\backslash I}z^{|I\backslash [k]|}
w^{|[n;m]\backslash I|}$.
\end{description}
Note that the polynomial $U_{0,m}^{(0)}$ coincides with Umemura's
polynomial $T_m(z^2,w^2;a,b)$. The formula for generalized
Umemura polynomials stated below follows from the Cauchy
identity, and was used by
 J. P. van Diejen and A. N. Kirillov [DK]
in their study of $q$-spherical functions.
\begin{Lemma}\label{l3.1}
The generalized Umemura polynomials $U_{n,m}^{(k)}(a,b;z,w)$
admit the following determinantal expression
\begin{eqnarray*}
U_{n,m}^{(k)}(a,b;z,w)=
\det \left|
a_i w^i \prod_{s \in [k]}\left(
\frac{i+s}{i-s}
\right) \delta_{i,j}
+\frac{2i}{i+j}\cc({i})
\prod_{s \in [n,m], s \neq i}
\left|
\frac{i+s}{i-s}
\right|
b_i z^i
\right|_{i,j\in [n;m] \backslash [k] },
\end{eqnarray*}
where $\cc(i)=i \ \ \mbox{if $i \leq n$, and }\ \  (i-n)/2
\ \ \mbox{if $i >n$}$.
\end{Lemma}
In the particular case $k=0$, $n=0$ this formula gives a determinantal
representation for Umemura's polynomials and has many applications.

\section{Main result}
\setcounter{equation}{0}
Let us introduce notation $U_{n,m}:=U_{n,m}^{(0)}(z,w;a,b)$.
The main result of our paper describes a
recurrence relation between polynomials
$U_{n,m}$.
\begin{Theorem}\label{t4.1}
\begin{eqnarray}
U_{n,m-1}U_{n,m+1}&=&
\left( -\bar{a}_{n+2m+2}z^2+\bar{b}_{n+2m+2}w^2\right)U_{n,m}^2+
8z^2w^2D_x^2 U_{n,m} \circ U_{n,m} \nonumber \\
\smallskip \nonumber\\
&&
-\frac{4}{(n+2m+1)^2}ab(a-b)z^2w^2\left( U_{n,m}^{(1)}\right)^2,
\end{eqnarray}
where for any two functions $f=f(x)$ and $g=g(x)$
\begin{eqnarray*}
D_x^2f \circ g=f''g-2f'g'+fg''
\end{eqnarray*}
denotes the second Hirota derivative, and $'=\frac{d}{dx}$;
here variables $z$,$w$ and $x$ are connected by the relations
$z=\frac{1}{2}\left( e^x+e^{-x}-2\right)^{1/2}$,
$w=\frac{1}{2}\left( e^x+e^{-x}+2\right)^{1/2}$.
\end{Theorem}
Below we give a sketch of our proof of Theorem~\ref{t4.1}.
Detailed exposition will appear elsewhere. The main step of the
proof is to establish the following algebraic identity which
appears to have an independent interest.
\begin{Lemma}\label{l4.2}
For any two subsets $I$, $J$ of the set $[n;m]$, we have
\begin{eqnarray*}
&&\prod_{ \lambda \in I}
\left(
\frac{ x + 2+\lambda}{ x +2- \lambda}
\right)
\prod_{ \lambda \in J}
\left(
\frac{ x - \lambda}{x  + \lambda}
\right)+
\prod_{ \lambda \in I}
\left(
\frac{ x - \lambda }{ x + \lambda }
\right)
\prod_{ \lambda \in J}
\left(
\frac{ x +2+ \lambda }{ x +2- \lambda }
\right)
\\
\smallskip \nonumber\\
&&\ \ \ = 2 + \sum_{ \lambda \in I \cup J}
\frac{b_{ \lambda }^{ I , J }}{( x +2- \lambda )( x + \lambda )},\nonumber
\end{eqnarray*}
where $b_{\lambda}^{I,J}$ are some constants (depending on $\lambda$, $I$
and $J$) which may be computed explicitly.
\end{Lemma}
One can prove this lemma by using the residue theorem.
\begin{Lemma}\label{l4.3}
For any two subsets $I$, $J$ of the set $[n;m]$, we have
\begin{eqnarray*}
\sum_{ \lambda \in I \cup J }
b_{ \lambda }^{ I , J } = 4( |I| - |J| )^2 - 4 (|I| + |J|).
\end{eqnarray*}
\end{Lemma}
This lemma follows from Lemma~\ref{l4.2}.
\begin{Lemma}\label{l4.4}
For an element $ \lambda \in I \cap J $, we have
$ b_{ \lambda } ^{ I , J }=0 $ if and only if $ \lambda - 2 \in I \cap J $.
For an element  $ \lambda \in I \backslash  (I \cap J) $, we have
$b_{ \lambda}^{I,J}=0 $ if and only if $ \lambda-2 \in J $.
\end{Lemma}
This lemma follows from Lemma~\ref{l4.2} by direct calculation.
\begin{flushright}
$\square$
\end{flushright}

\begin{description}
\item[Remarks 1.]If $n=0$, then $U_{0,m}=U_{m+1}(z^2,w^2;a,b)$
coincides with the Umemura polynomial, and $U_{0,m}^{(1)}=0$. In
this case the recurrence relation (4.1) has been used by M. Taneda
in his proof of Okada--Noumi--Okamoto--Umemura's Conjecture (2.6).
\item[2.]
Note that $U_{0,m}=U_{2,m-1}^{(1)}/(2m+1)$, and more generally\\
$U_{k,m}^{(k)}=U_{k+2,m-1}^{(k+1)}(2k+1)!!(2m-1)!!/(2k+2m+1)!!,$
where $(2n+1)!!=1 \cdot 3 \cdot 5 \cdots (2n+1)$.
\item[3.]
"Unwanted term" in (4.1) which contains $\left( U_{n,m}^{(1)} \right)^2$
vanishes if either $a=0$, or $b=0$, or $a=b$.
\end{description}
In the case $a=b$ and $k=0$ the expression $e_I^{(n,m,k)}(z,w)$
doesn't depend on a subset $I \subset [n;m]$ and is equal to
$a_{[n;m]} z^{|I|}w^{|[n;m] \backslash I|}$. Hence,
in this case
\begin{eqnarray}
U_{n;m}^{(0)}(z,w;a,a)&=&
a_{[n;m]} \sum_{I \subset [n;m]} \dd_{n,m}(I) (-1)^{\cc(I)}
z^{|I|}w^{|[n;m] \backslash I|} \nonumber \\
\smallskip \nonumber \\
&=&a_{[n;m]}(z+w)^{\tiny  \left(
\begin{array}{c} n+m+1 \\2
\end{array}\right)}
(z-w)^{\tiny \left( \begin{array}{c} m+1 \\2
\end{array}\right) }.
\end{eqnarray}
Recall that $\bar{a}_i=a+(i-1)^2$, $a_{2i}=\bar{a}_2 \bar{a}_4
\cdots \bar{a}_{2i}$, $a_{2i+1}=\bar{a}_1 \bar{a}_3 \cdots
\bar{a}_{2i+1}$ and $a_{[n,m]}=\prod_{ i \in [n;m] } a_i$. The
last equality in (4.2) has been proved for the first time by
J.P.~van Diejen and A.N.~Ki\-rillov [DK]. On the other hand, we
can show that polynomials
\begin{eqnarray*}
X_{n,m}(z,w;a)=a_{[n;m]}(z+w)^{\tiny  \left(
\begin{array}{c} n+m+1  \\2
\end{array}\right)}
(z-w)^{\tiny \left( \begin{array}{c} m+1 \\2
\end{array}\right) }
\end{eqnarray*}
also satisfy the recurrence relation (4.1) and coincide with
polynomials
$U_{n,m}^{(0)}(z,w;a,a)$ if $m=0$. From this observation we can deduce the
equality
$X_{n,m}(z,w;a)=U_{n,m}^{(0)}(z,w;a,a)$, which is equivalent to the main
identity from [DK].
Another case when "unwanted term" in (4.1) vanishes is the case when either
$a=0$, or $b=0$.
In this case we have
\begin{Corollary}\label{c4.5}
Assume that $b=0$, then polynomial $U_{n,m}(z,w;a,0)$
defines a solution to the equation Painlev\'e $VI$.
\end{Corollary}
Finally, we are going to compare polynomials $U_{n,m}(z,w;a,0)$ and
$U_m(z,w;\alpha,\beta)$.
For this goal, let us consider functions
\begin{eqnarray*}
h_0:=h_0(t)=\left\{ b_1^2\left( \sqrt{t}-\sqrt{t-1} \right)^2
+b_2^2\left( \sqrt{t}+\sqrt{t+1}\right)^2
\right\}/4,
\end{eqnarray*}
and
\begin{eqnarray*}
h_{n,m}:=h_{n,m}(b_1,b_2)=t(t-1) \log \left(U_{n,m} \right)'-h_0.
\end{eqnarray*}
\begin{Proposition}\label{p4.6}
\begin{description}
\item[(i)]  $h_{0,m}$ satisfies the Painlev\'e--Okamoto
equation $E_{VI}(b_1,b_2,m+1/2,0)$;
\item[(ii)] $h_{1,m}=-(2t-1)(m+1)^2/2$ satisfies
the equation $E_{VI}(0,m+1,b_3,b_4)$;
\item[(iii)] $h_{n,m}(0,b_2)$ satisfies the equation
$E_{VI}(0,b_2,\frac{n}{2},\frac{n+2m+1}{2})$.
\end{description}
\end{Proposition}
Proposition~\ref{p4.6} follows from Lemma~\ref{l4.7} and
Lemma~\ref{l4.8} below. Let us define
$U_{n,m}(b_1,b_2):=U_{n,m}(z,w;-4b_1^2,-4b_2^2)$, then
\begin{Lemma}\label{l4.7}
\begin{eqnarray*}
U_{n,m}(0,b_2)
=\left\{
\begin{array}{lcrc}
\displaystyle b_{[n;m]_\odd} w^{\left(\frac{n}{2}\right)^2}
U_{0,m+\frac{n}{2}}\left(\frac{n}{2},b_2\right), & \mbox{if $n$
is even},
&(4.3)\\
\smallskip \\
\displaystyle b_{[n;m]_\odd} w^{\left(\frac{n+2m+1}{2}\right)^2}
U_{0,\frac{n-1}{2}}\left(m+\frac{n+1}{2},b_2\right), & \mbox{if
$n$ is odd},  &(4.4)
\end{array}
\right.
\end{eqnarray*}
where $[n;m]_\odd=\left\{ i \in [n;m] \ \ |\ \  \mbox{$i$ is odd }\right\}$.
\end{Lemma}
From Lemma~\ref{l4.7} we can deduce the following
\begin{Lemma}\label{l4.8}
\begin{eqnarray*}
h_{n,m}(0,b_2)
=\left\{
\begin{array}{lc}
\displaystyle
h_{0,m+\frac{n}{2}}\left(\frac{n}{2},b_2\right),& \mbox{if $n$ is even,}\\
\smallskip \\
\displaystyle
h_{0,\frac{n-1}{2}}\left(m+\frac{n+1}{2},b_2\right),
&\mbox{if $n$ is odd.}
\end{array}
\right.
\end{eqnarray*}
\end{Lemma}
It follows from Lemma~\ref{l4.7}, (4.4) and Theorem~\ref{t4.1},
that Umemura's polynomials $U_{m}(b_1,b_2)$ satisfy a new
recurrence relation with respect to the first argument $b_1$.
\begin{Theorem}\label{t4.9}
\begin{eqnarray*}
\lefteqn{
U_m(b_1-1,b_2)U_m(b_1+1,b_2)(b_1^2-b_2^2)
}\\ \smallskip \\
&=&
(b_1^2-b_2^2)U_m^2(b_1,b_2)+2z^2D_x^2 U_m(b_1,b_2) \circ U_m(b_1,b_2).
\end{eqnarray*}
Recall that $D_x^2$ denotes the second Hirota derivative.
\end{Theorem}

\section{Conjecture}
\setcounter{equation}{0}
We define $q_m:=q_m(t)$ by
\begin{eqnarray*}
q_m-t&=&4U_m^2\left\{
\left(m+\frac{1}{2}\right)t(t-1)\frac{d}{dt} \log U_{m+1}
-
\left(m+\frac{3}{2}\right)t(t-1)\frac{d}{dt} \log U_{m}\right.\\
\smallskip \nonumber\\
&&
\left.
-\frac{1}{2} b_1 b_2
+\frac{1}{4}\left(b_1^2\frac{z}{w}+b_2^2\frac{w}{z}\right)
\right\} /\left(U_{m+1}U_{m-1}-(2m+1)^2U_m^2\right).
\end{eqnarray*}

One can check that $q_m$ is a solution to both equations
$P_{VI}(b_1,b_2,m+\frac{1}{2},0)$ and\break
$P_{VI}(b_1,b_2,0,m+\frac{1}{2})$. It follows from Okamoto's
theory [OI] that the function
\begin{eqnarray}
\bar{h}_{1,m}=t(t-1)\frac{d}{dt} \log U_{m+1} -
\frac{1}{4}\left(b_1^2\frac{z}{w}+b_2^2\frac{w}{z}\right)
+\left( m+\frac{1}{2}\right) q_m -\frac{1}{2}\left( m+ \frac{1}{2}\right)
\end{eqnarray}
is also a solution to $E_{VI}(b_1,b_2,n+\frac{1}{2},1)$.
Based on the latter expression for the function $\bar{h}_{1,m}$, and
using Lemmas 5 and 6, we come to the following
\begin{Conjecture}\label{c5.1}
If $b_1=0$, then we have
\begin{eqnarray*}
U_{m+1}U_{m-1}-(2m+1)^2U_m^2=\frac{1}{4b_2^2}U_{2,m-1}^2,
\end{eqnarray*}
where $U_m:=U_m(0,b_2)$ is a special case of Umemura's polynomial,
and $U_{2,m-1}=U_{2,m-1}(0,b_2)$.
\end{Conjecture}

\begin{flushleft}
\vskip 0.5cm
Anatol N. KIRILLOV \\
Graduate School of Mathematics,
Nagoya University\\

Chikusa-ku, Nagoya,464-8602, Japan\\

and\\

Steklov Mathematical Institute \\

Fontanka 27, St Petersburg 191011, Russia

\begin{verbatim}
E-mail address: kirillov@math.nagoya-u.ac.jp
\end{verbatim}
\vskip 0.3cm
Makoto, TANEDA\\
Hida 4-3-87, Kumamoto, 861-5514, Japan.\\
\begin{verbatim}
E-mail address: tane@rc4.so-net.ne.jp
\end{verbatim}
\end{flushleft}

\end{document}